\documentclass[12pt]{article}
\usepackage{amsxtra,amssymb,amsthm,amsmath,latexsym}

\def\R{{\mathbb R}}
\def\q{\quad}

\def\oH{\buildrel\circ\over H}
\def\oH1{\buildrel\circ\over H\kern-.02in{}^1}

\def\f{\frac}
\def\s{\sigma}
\def\b{\beta}
\def\ep{\epsilon}

\def\d{{\delta}}

\begin{document}
 A.G. Ramm, Global convergence for ill-posed equations with monotone 
operators: the dynamical systems method, J.Phys A, 36, (2003), L249-L254

\title{
Global convergence for ill-posed equations with monotone operators: 
the dynamical systems method
\footnote{Math subject classification: 34R30,  35R25, 35R30, 
37C35, 37L05, 37N30, 47A52,  
47J06, 65M30, 65N21 }
}

\author{ A.G. Ramm\\
Mathematics Department, Kansas State University, \\
Manhattan, KS 66506-2602, USA\\
Email: ramm@math.ksu.edu
}

\date{}

\maketitle\thispagestyle{empty}

{\it  Consider an operator equation $F(u)=0$ in a real Hilbert space.
 Let us call this equation ill-posed if the 
operator
 $F'(u)$ is not boundedly invertible, and well-posed otherwise.
 If $F$ is monotone $C^2_{loc}(H)$ operator, then we construct a Cauchy 
problem, which has the following properties:
 1) it has a global solution for an arbitrary initial data,
 2) this solution tends to a limit as time tends to infinity,
 3) the limit is the minimum norm solution to the equation $F(u)=0$.

Example of applications to linear ill-posed operator equation is given.
}

\section{Introduction}
Many physical problems can be formulated as operator equations.
In this paper a general convergence theorem is proved
for solving operator equations with monotone operators.
Consider an operator equation
$$
F(u):=\mathcal{B}(u)-f=0, \quad f\in H,
\eqno{(1.1)}
$$
where $\mathcal{B}$ is a monotone, nonlinear, $C^2_{loc}$ operator in a 
real Hilbert 
space $H$, i.e, $\sup_{u\in B(u_0,R)}||F^{(j)}(u)||\leq M_j(R):=M_j,\,j=0, 
1, 2,$ where $R>0$ is arbitrary,
$B(u_0,R):=\{u: ||u-u_0||\leq R \}$, $F^{(j)}(u)$ is the Fr\'echet
derivative. 
Let $N:=\{z: F(z)=0\}$. It is known
that $N$ is convex and closed under our assumptions.
Assume that $N$ is not empty. Then it contains the unique minimum norm 
element $y$: $F(y)=0, \quad ||y||\leq ||z||,\,\, \forall z\in N.$
These assumptions hold throughout and are not repeated.

Let $\dot u$ denote derivative with respect to time. Consider the 
dynamical system ( the Cauchy problem ):
$$
\dot u=\Phi(t,u), \,\,\, u(0)=u_0;\,\, \Phi:=-A_\ep^{-1}[F+\ep u], 
\eqno{(1.2)}
$$
where $A_\ep:=A+\ep I$, $A:=F'(u)$, $I$ is the identity operator,
and $\ep=\ep(t)>0$ is a continuously differentiable, monotone, decaying
to zero as $t\to \infty$, function on $[0, \infty)$. 
Specifically, we will use $\ep=c_1(c_0+t)^{-b},$ where 
$c_1,c_0$ and $b$ are positive constants, $b\in (0,1)$,
and assume throughout (without repeating), that $|\dot \ep|\ep^{-1}\leq 
0.25$. Note that
$\Phi(t,u)$ is locally Lipschitz with respect to $u\in H$ and 
continuous with respect to $t\geq 0$ under our assumptions.
Thus problem (1.2) has a unique local solution.
We want to solve equation (1.1) by  
solving (1.2), and proving that for any initial $u_0$ the following 
three results hold:
$$
\exists u(t) \forall t>0; \quad \exists u(\infty):=\lim_{t\to \infty}u(t); 
\quad F(u(\infty))=0.
\eqno{(1.3)}
$$
Moreover, we prove that the solution $u\in B(u_0,R)\,\, \forall t\geq 0,$
where $R:=3r$, and $r:=||y||+||u_0||$.

Problem (1.1) with noisy data
$f_\delta$, $||f_\delta -f||\leq \delta$, given in place of $f$, generates
the problem:
$$
\dot u_\d=\Phi_\d(t,u_\d), \,\,\, u_\d(0)=u_0,
\eqno{(1.4)}
$$
The solution $u_\d$ to (1.4), calculated at
a suitable stopping time $t=t_\d$, converges to $y$:
$$
\lim_{\d \to 0}||u_\d(t_\d)-y||=0.
\eqno{(1.5)}
$$
 The choice of $t_\d$ with this property is called the stopping rule.
One has usually $\lim_{\d \to 0}t_\d=\infty$.

We do not restrict the growth of nonlinearity at infinity and do not 
assume that the initial approximation $u_0$ is close to the solution $y$
in any sense. Usually (e.g., see [2]) convergence theorems for
Newton-type methods for solving nonlinear equation (1.1) 
have the assumption that the initial data $u_0$ is close to $y$.
We obtain a global convergence result for a 
continuous regularized Newton-type method 
(1.2). This result is stated in Theorem 1, and proved in Section 2. 

{\bf Theorem 1.} {\it For any choice of $u_0$
problem (1.2) has a global solution, this solution stays in the ball
$B(u_0, R)$, and (1.3) holds.  
If $u_\d(t)$ solves (1.4),
then there is a $t_\d$ such that $\lim_{\d \to 0}||u_\d(t_\d)-y||=0$.
}

The proof uses essentially the following result which is obtained in [1].

{\bf Theorem 2.} {\it Let $\gamma(t), \s(t), \b(t)\in C[t_0,\infty)$ for
some real number $t_0$.
If there exists a positive  function $\mu(t)\in C^1[t_0,\infty)$ such that
$$
0\le\s(t)\le\f{\mu(t)}{2}[\gamma(t)-\frac{\dot{\mu}(t)}{\mu(t)}],\quad
\b(t)\le\frac{1}{2\mu(t)}[\gamma(t)-\frac{\dot{\mu}(t)}{\mu(t)}],\quad
g_0\mu(t_0)<1,
\eqno{(1.6)}
$$
where $g_0$ is the initial condition in (1.7), then a nonnegative
solution $g$ to the following differential inequality:
$$
\dot{g}(t)\le -\gamma(t) g(t)+\s(t)g^2(t)+\b(t),\quad
g(t_0)=g_0,
\eqno{(1.7)}
$$
satisfies the estimate:
$$
0\leq g(t)\,\le\, \f{1-\nu(t)}{\mu(t)}\,<\, \f{1}{\mu(t)},
\eqno{(1.8)}
$$
for all $t\in [t_0,\infty)$, where
$$
0<\nu(t)=\left(\f{1}{1-\mu(t_0)g(t_0)}+
\f{1}{2}\int_{t_0}^t\left(\gamma(s)-\f{\dot{\mu}(s)}{\mu(s)}\right)ds
\right)^{-1}.
\eqno{(1.9)}
$$
} 
There are several novel features in this result. First, differential
equation, which one gets from (1.7) by replacing the inequality sign
by the equality sign, is a Riccati equation, whose solution may blow up in
a finite time, in general. Conditions (1.6) guarantee the
global existence of the solution to this Riccati equation with the initial
condition (1.7). Secondly, this Riccati differential equation
cannot be integrated analytically by separation of variables, in general.
Thirdly, the coefficient $\sigma(t)$ may grow to infinity as
$t\to \infty$, so that the quadratic term does not necessarily has a
small coefficient, or the coefficient smaller than $\gamma(t)$.
Without loss of generality one may assume $\beta(t)\geq 0$ in Theorem 2.
 In [4] and [5] one finds a description and applications of DSM
(dynamical systems method) and some 
remarks about discrepancy principle, which are useful in treating
problems with noisy data. Many physical problems can be formulated as 
operator equations with monotone operators. We mention the theory of 
passive networks (see [7] and [8], Chapter 3) as just one of many 
examples.

 \section{Proof of Theorem 1}

Let us sketch the proof.
Denote $w:=u-V$, $||w||:=g$, $v:=||V-y||$. Clearly
$||u(t)-y||\leq g+v$. We will prove that $\lim_{t\to \infty} g=\lim_{t\to 
\infty} v=0$.

Let $V$ solve the equation
$$
F(V)+\ep(t) V=0.
\eqno{(2.1)}
$$
Under our assumptions on $F$, it is known that (2.1) has a unique 
solution for every $t>0$, and $\lim_{t \to \infty}||V(t)-y||=0.$ 
One can prove that $\sup_{t\geq 0}||V||\leq ||y||$,
$V$ is differentiable, and $||\dot V(t)||\leq 
||y|| |\dot \ep(t)|/\ep(t)$. 
We will show that the global solution $u$ to (1.2)
does exist, and $\lim_{t \to \infty}||u(t)-V(t)||=0$. 
This is done by deriving a differential 
inequality for $w$, and by applying Theorem 2 to $g=||w||$.
Since $||u(t)-y||\leq g+v$, it then follows that
(1.3) holds. We also check that $u(t)\in B(u_0, R)$, where 
$R:=3(||y||+||u_0||)$, for any choice of $u_0$ and a suitable choice of 
$\ep=\ep(t)$.

Let us derive the differential inequality for $w$. One has
$$\dot w=-\dot V -A_{\ep(t)}^{-1}(u)\bigl[F(u(t))-F(V(t))+\ep(t)w],
\eqno{(2.2)}
$$
and $F(u)-F(V)=Aw+K$, where $||K||\leq M_2g^2/2$, $g:=||w||$.
Multiply (2.2) by $w$, use monotonicity of $F$, i.e.,
the property $A\geq 0$,
and the estimate $||\dot V ||\leq ||y|| |\dot \ep|/\ep$, and get:
$$
\dot g\leq -g +\frac {0.5 M g^2}{\ep} +||y||\frac{|\dot \ep|}{\ep}, 
\eqno{(2.3)}
$$
where $M:=M_2$. 
Inequality (2.3) is of the type (1.7): $\gamma=1$, $\s=0.5 M/\ep$,
$\b=||y||\frac{|\dot \ep|}{\ep}$. Choose $\mu(t)=2M/\ep(t)$.
 Clearly $\mu\to \infty$ as $t\to \infty$.
Let us check three conditions (1.6). One has 
$\frac{\dot{\mu}(t)}{\mu(t)}=|\dot \ep|/\ep$. Take $\ep=
c_1(c_0+t)^{-b}$, where $c_j>0$ are constants, $0<b<1$,
and choose these constants so that $|\dot \ep|/\ep <1/2,$
e.g., $\frac b{c_0}=\frac 1 4$.  
Then the first condition (1.6) is satisfied.
The second condition (1.6)
holds if (*) $8M||y|| |\dot \ep|\ep^{-2}\leq 1$. 
One has $\ep(0)=c_1c_0^{-b}$. Choose $\ep(0)=4Mr$.
Then
$ |\dot \ep|\ep^{-2}=bc_1^{-1}(c_0+t)^{b-1}\leq 
bc_0^{-1}c_1^{-1}c_0^b=\frac 1 {4\ep(0)}=\frac 1 {16Mr}$. 
Thus, the second condition (1.6) holds.
The last
condition (1.6) holds because  
$2M||u_0-V_0||/\ep(0)\leq \frac {2Mr}{4Mr}=\frac 1 2<1$.

Thus, by Theorem 2,
$g=||w(t)||<\frac {\ep(t)}{2M}\to 0$ when $t\to \infty$, and
$||u(t)-u_0||\leq g+||V-u_0||\leq g(0)+r\leq 3r$. This estimate 
implies the global 
existence of the solution to (1.2), because if $u(t)$ would have a finite
maximal interval of existence, $[0,T)$, then $u(t)$ could not stay bounded 
when $t\to T$, and this contradicts the boundedness of $||u(t)||$,
which follows from our estimates: $||u(t)||\leq 4r$.
We have proved the first part of Theorem 4.2, namely properties (1.3).  
$\Box$

To derive a stopping rule we argue as follows. One has
$||u_\d(t)-y||\leq ||u_\d(t)-V(t)||+||V(t)-y||:=g_\d +v$.
We have proved that $\lim_{t\to \infty} v(t):=0$.
The rate of decay of $v(t)$ can be arbitrarily slow, in general.
Additional assumptions, e.g., the source-type ones, can be used to 
estimate the rate of decay of $v(t)$. One can derive differential 
inequality similar to (2.3)
 for $g_\d:=||u_\d(t)-V(t)||$, and estimate 
$g_\d$ using (1.8). The analog of (2.3) for $g_\d$ contains additional
term $\d/\ep$ on the right-hand side. If $16M\d \leq \ep^2,$ then
conditions (1.6) hold, and $g_\d<\frac {\ep(t)} {2M}$. 
Let $t_\d$ be the root of the equation $\ep^2(t)=16M\d$.
Then $\lim_{\d\to 0}t_\d=\infty$, and 
$\lim_{\d \to 0}||u_\d(t_\d)-y||=0$,
because $||u_\d(t_\d)-y||\leq v(t_\d)+g_\d$, 
$\lim_{t_\d\to \infty}g_\d(t_\d)=0$ and 
$\lim_{t_\d\to \infty} v(t_\d) =0$, but the convergence can be 
slow. See also [3] for the rate of convergence under source assumptions.
If the rate of decay of $v(t)$ is known, then one can
choose $t_\d$ 
as the minimizer of the problem, similar to (3.13), 
$v(t)+g_\d(t)=min$ , where the minimum is 
taken over $t>0$ for a fixed small $\d>0$. This would yield a quasioptimal 
stopping rule.
Theorem 1 is proved. $\Box$ 

\section{Example}

Let us give an example of applications of Theorem 1.
Consider a linear operator equation: 
$$
Au=f.
\eqno{(3.1)}$$
Let us denote by A) the folowing assumption:

{\bf Assumption} (A):  {\it $A$ is a linear, bounded operator in
$H$, defined on all of $H$,
the range $R(A)$ is not closed, so (3.1) is an ill-posed
problem, there is a $y$ such that $Ay=f$, $y\perp N$, where $N$ is the
null-space of $A$. }

Let $\mathcal{B}:=A^*A$,
$q:=\mathcal{B}y=A^*f$, $A^*$ is the adjoint of $A$.
Every solution to (3.1) solves
$$
\mathcal{B}u=q,
\eqno{(3.2)}$$
and, if $f=Ay$, then every  solution to (3.2) solves (3.1). Choose a
continuous function  $\ep(t)>0$, monotonically
decaying to zero on $\R_+$, as in Theorem 1.
If $\mathcal{B}$ is a linear operator, and $F(u):=\mathcal{B}u-q$, then
$F'(u)=\mathcal{B}$, and 
$\Phi:=-(\mathcal{B}+\ep)^{-1}[\mathcal{B}u-q+\ep u]=
-u+(\mathcal{B}+\ep(t))^{-1}q$. Therefore equation (1.2)
takes the form:
$$
\dot u=-u+(\mathcal{B}+\ep(t))^{-1}q,\quad u(0)=u_0.
\eqno{(3.3)}
$$
The operator $\mathcal{B}:=A^*A\geq 0$ is linear, monotone, and Theorem 1 
is applicable. Therefore conclusions (1.3) hold for the solution to
(3.3),  and, since equations (3.1) and (3.2) are equivalent
if (3.1) is solvable, one concludes that $u(\infty)=y$, where $y$ is the
unique minimal-norm solution to equation (3.1).
Moreover, if the data are noisy, so that $f_\d$ is given in place of 
$f$, and $||f-f_\d||\leq \d$, then Theorem 1 yields a stable solution
to the ill-posed problem (3.1).
Thus, Theorem 1 yields a method for solving arbitrary linear
ill-posed problems with bounded linear operator $A$.
This method works well numerically.

\section{Appendix}
 For convenience of the reader and for completeness of the presentation we 
include a proof of Theorem 2 which is borrowed from [1].

We start with the well-known lemma (see, e.g. [6]):

{\bf Lemma.}
{\it Let $f(t,w)$, $g(t,u)$ be continuous on region $[0,T)\times D$ 
($D\subset
R$, $T\le\infty$) and $f(t,w)\leq g(t,u)$ if $w\leq u$,  $t\in (0,T)$,
$w,u\in D$. Assume that  $g(t,u)$ is such that the Cauchy problem
$$
\dot u = g(t,u), \q u(0)=u_0,\q u_0\in D
$$
has a unique solution. If
$$
\dot w\leq f(t,w), \q w(0)=w_0\leq u_0,\q w_0\in D,
$$
then $u(t)\geq w(t)$ for all $t$ for which $u(t)$ and $w(t)$ are
defined.
}

Let us now prove Theorem 2.

{\bf Proof of Theorem 2}.

Let $g$ be the function from (1.7). Define the new function $w$ by the 
formula:
$$
w(t):=g(t)e^{\int_{t_0}^t \gamma(s)ds}.
$$
Then 
$$
\dot{w}(t)\le a(t)w^2(t)+b(t),\q w(t_0)=g(t_0),
$$
where
$$
a(t)=\s(t)e^{-\int_{t_0}^t \gamma(s)ds},\q b(t)=\b(t)e^{\int_{t_0}^t 
\gamma(s)ds}
$$
Consider the equation:
$$
\dot{u}(t)=\f{\dot{f}(t)}{G(t)}u^2(t)-\f{\dot{G}(t)}{f(t)}.
\eqno{(4.1)}
$$
One can check by a direct calculation that the the solution to
this equation is given by the following formula (see [6]):
$$
u(t)=-\f{G(t)}{f(t)}+\left[f^2(t)\left(C-
\int_{t_0}^t\f{\dot{f}(s)}{G(s)f^2(s)}ds\right)\right]^{-1}, 
\eqno{(4.2)}
$$
where $C$ is a constant. If $u(0)=u_0$, then $C=\frac 1 
{u_0f^2(0)+G(0)f(0)}$.
Define $f$ and $G$ as follows:
$$
f(t):=\mu^{\f{1}{2}}(t)e^{-\f{1}{2}\int_{t_0}^t \gamma(s)ds},\q
G(t):=-\mu^{-\f{1}{2}}(t)e^{\f{1}{2}\int_{t_0}^t \gamma(s)ds},
$$
and consider the Cauchy problem for equation
(4.1) with the initial condition
$u(t_0)=g(t_0)$. Then $C$ in (4.2) can be calculated:
$$C=\f{1}{\mu(t_0)g(t_0)-1}.$$
From (1.6) one gets
$$
a(t)\le \f{\dot{f}(t)}{G(t)},\q b(t)\le -\f{\dot{G}(t)}{f(t)}.
$$
Since $fG=-1$ one has:
$$
\int_{t_0}^t\f{\dot{f}(s)}{G(s)f^2(s)}ds=-\int_{t_0}^t\f{\dot{f}(s)}{f(s)}ds=
\f{1}{2}\int_{t_0}^t\left(\gamma(s)-\f{\dot{\mu}(s)}{\mu(s)}\right)ds.
$$
Thus
$$
u(t)=\f{e^{\int_{t_0}^t \gamma(s)ds}}{\mu(t)}\left[1-\left(\f{1}
{1-\mu(t_0)g(t_0)}+
\f{1}{2}\int_{t_0}^t\left(\gamma(s)-\f{\dot{\mu}(s)}{\mu(s)}\right)ds
\right)^{-1}\right].
\eqno{(4.3)}
$$
It follows from conditions (1.6) and from the second inequality in
(1.6) that the solution to
problem (4.1) exists for all $t\in[0,\infty)$ and the following 
inequality
holds with $\nu(t)$ defined by (1.9):
$$
1\, >\, 1-\nu(t)\ge\mu(t_0)g(t_0).
$$
From Lemma and from formula (4.3) one gets:
$$
g(t)e^{\int_{t_0}^t \gamma(s)ds}:=w(t)\le u(t)=\f{1-\nu(t)}{\mu(t)}
e^{\int_{t_0}^t \gamma(s)ds}<\f{1}{\mu(t)}e^{\int_{t_0}^t\gamma(s)ds},
$$
and thus estimate (1.8) is proved.
Theorem 2 is proved. $\Box$

\end{document}